# Smoothing Before Estimating Uncertainty, Scaling, and Intermittency: Application to Short Heart Rate Signals

David R. Bickel[1]

*Office of Biostatistics and Bioinformatics, Medical College of Georgia*
*1120 Fifteenth St., AE-3037, Augusta, Georgia 30912-4900*

August 23, 2002

**Abstract.** Three aspects of time series are uncertainty (dispersion at a given time scale), scaling (time-scale dependence), and intermittency (inclination to change dynamics). Simple measures of dispersion are the mean absolute deviation and the standard deviation; scaling exponents describe how dispersions change with the time scale. Intermittency has been defined as a difference between two scaling exponents. After taking a moving average, these measures give descriptive information, even for short heart rate records. For this data, dispersion and intermittency perform better than scaling exponents.

[1] Electronic addresses: dbickel@mail.mcg.edu, bickel@prueba.info; URL: www.davidbickel.com

*Introduction.* Healthy and heart-failure subjects have been distinguished by various measures of heart rate variability (HRV), including measures of dispersion, scaling exponents, intermittency, and multifractality. These four classes of estimators are related: dispersion is a scale-dependent measure of uncertainty such as the standard deviation, scaling exponents specify how dispersion depends on the time scale, and intermittency and multifractality can be computed from multiple scaling exponents. One of the first estimators of HRV is the standard deviation of the heart rate; it is now well-known that a low standard deviation is a sign of poor health. West and Goldberger [1] hypothesized that pathology would lead to a narrowing of the power spectrum of the heart rate; this would imply differences in scaling exponents between sick and healthy subjects. Such differences have been confirmed [2], but certain measures of dispersion appear to discriminate better between the two classes of subjects, especially for shorter time series [3]. For time series that are several hours long, measures of intermittency [4] and multifractality [2] also distinguish subjects with congestive heart failure from healthier subjects. However, the measurement of continuous heart rate data for such a long time is often impractical clinically since it either requires that the patients are bedridden and attached to hospital monitors, or that they wear expensive ambulatory monitors, as in Ref. 2. Many research protocols also limit the duration of monitoring: the heart rate data of Ref. [5] and of this Letter could only be collected while the nurse researchers were physically present. In addition, heart rate records that are freely available tend to be short; the public records from www.physionet.org that were analyzed in Ref. [6] are each only 30 minutes long. In this Letter, methods for estimating the dispersion, scaling exponents, and intermittency that apply to short physiological records as well as long ones are proposed and evaluated; these methods can also be used to describe other fractal time series. The method of multifractality quantification of Ref. 2 is not treated herein since it is computationally less efficient, requiring the estimation of many scaling exponents, and since it is probably inappropriate for short time series.

The $q$th power deviation of a time series $(x_1, x_2, x_3,...)$ is defined as

$$\sigma_k(q) = \langle |x_{k+i} - x_i|^q \rangle^{1/q}, \qquad (1)$$

where $k$ is the positive integer scale and $q$ is the positive real order. The angular brackets denote the arithmetic mean. $\sigma_k(q)$ meets all three of the requirements of a measure of dispersion listed in Ref. 7. $\sigma_k^q(q)$ is the standard structure function; wavelet-based measures of dispersion have been used to handle time series without stationary increments [2,3], but the time series considered herein are short enough that $x_{i+1} - x_i$ can be assumed to be stationary. Further, wavelet-based estimates have higher variances [8], which aggravate the uncertainty already inherent in estimates from short time series, so the simpler estimator based on the standard structure function is used in this Letter. Even so, the proposed preprocessing techniques can also be applied before computing wavelet-based estimates of dispersion, scaling, and intermittency when such estimates are appropriate. The generalized Hurst exponent, $H(q)$, is the scaling exponent defined by $\sigma_k(q) \propto k^{H(q)}$ and is estimated as the slope of the line that best fits a plot of $\log \sigma_k(q)$ versus $\log k$ for some values of $k$; $H(q)$ was used in Ref. [9]. In general, $0 \le H(q) \le 1$ and, for a random walk that is the integral of Gaussian, uncorrelated noise, $\forall_{q>0} H(q) = 1/2$. Since $H(q)$ is also independent of $q$ for other non-intermittent processes, the intermittency of the time series can be quantified by a normalized difference in generalized Hurst exponents,

$$\chi(q_1, q_2) = -q_1 q_2 \frac{H(q_2) - H(q_1)}{q_2 - q_1}, \qquad (2)$$

called the *finite-difference intermittency*; the normalization factor ensures that $0 \le \chi(q_1, q_2) \le 1$ when time averages are used [4]. $\chi(1,q)$ is related to the generalized dimensions [4] and $\chi(1,2)$ is a measure of the intermittency of stationary point processes, with applications to DNA composition and to muscle activity [10].

*Smoothing.* High-frequency noise is often removed from time series before the estimation of parameters of interest; for example, Wang [11] used wavelet shrinkage to remove noise from a time series before estimating scaling exponents. A simpler and more common technique for reducing noise is the *moving average*, which smooths the data. For an odd, positive window size $m$, the $m$-point moving average of a time series replaces each value of the series with the mean of that value, the previous $(m-1)/2$ values, and the next $(m-1)/2$ values; the first $(m-1)/2$ and last $(m-1)/2$ values of the original time series are discarded rather than replaced. In order to preserve scaling properties, I propose applying the moving average to the nonstationary signal $(x_1, x_2, x_3,...)$, so that high-frequency fluctuations in the stationary increments of the signals are not necessarily removed. In addition to reducing high-frequency Gaussian noise in the nonstationary signal, the moving average reduces the effect of isolated outliers, values that are far from other values in the time series.

The signal can also be smoothed using the median instead of the mean; this *moving median* is much more resistant to outliers. An advantage of the moving median for some applications is that it resists the influence of isolated values while preserving sustained increases and decreases in a time series, e.g., the 3-point moving median of the time series $(0,0,0,1,0,0,1,1,0,1,1)$ is $(0,0,0,0,0,1,1,1,1)$, whereas the 3-point moving average is $(0,1/3,1/3,1/3,1/3,2/3,2/3,3,2/3,2/3)$. The resistance of an estimator to outliers can be quantified by the *breakdown point*, the smallest fraction of values replaced by arbitrary values that can make the estimator unbounded [12]. For example, the mean of $m$ values has a breakdown point of $1/m$ since, given any bound, the mean will be higher than that bound if one of the values is replaced by a sufficiently high value. The median, on the other hand, has an asymptotic breakdown point of 50%: at least half of the values must be replaced by outliers to move the median past any bound. (Other estimators of location, including certain estimators of the mode [13,14], also have breakdown points of 50%, but these estimators are not as easily computed as the median.) The high breakdown point of the median implies that the moving median is highly resistant to outliers [15], whereas the moving average is very sensitive to outliers.

That the moving median is more resistant to outliers than the moving average does not imply the superiority of the former for observed data since what appear to be outliers may represent important aspects of the dynamics of the process. In addition, for *m* observations drawn independently from a normal distribution of standard deviation $\sigma$, the standard deviation of estimates of the mean is asymptotically $\sigma/\sqrt{m}$, lower than that of the estimates of the median, $\sigma\sqrt{\pi}/\sqrt{2m} \approx 1.253\sigma/\sqrt{m}$, as $m \to \infty$ [16]. Thus, the *m*-point moving average is more effective than the *m*-point moving median in removing high-frequency white noise, reducing its standard deviation by a factor of approximately $\sqrt{m}$, rather than only $\sqrt{m}/1.253$. For most nonlinear systems, data cannot be neatly decomposed into signal, noise, and outliers, and the applicability of one smoothing function over another depends on the particular data set and the properties to be estimated. Consequently, I use heart rate data to compare the performance of both methods of smoothing in estimating dispersion, scaling exponents, and intermittency.

The moving average and moving median techniques obviate more complicated outlier-detection steps that are heart-rate-specific, such as those used in Ref. [2], and have been tested and optimized, as described below. Scaling in data can also be quantified in the presence of outliers by generating a sequence of symbols from continuous data [17], but the loss of information involved in that method precludes the estimation of intermittency.

*Application to heart rate data.* The proposed methods of fractal time-series analysis were applied to the times between successive heart beats, called interbeat intervals (IBIs), for 23 "stroke" subjects who have had a stroke and who have ischemic or hemorrhagic brain disease and for 23 "cardiac" subjects who have angina or ischemic heart disease, while they were at rest for about 8 minutes (a mean of 297 IBIs per subject). Originally, 24 stroke and 25 cardiac Japanese participants over the age of 60, with roughly equal numbers of males and females, had been recruited for measurement, but 3 were eliminated due to equipment failure or the possibility that a subject was not at rest (calling into question the assumption of stationary increments). To correct for the effect of different mean heart rates on dispersion estimates, each IBI was divided by the mean of all IBIs for that subject before applying the above analyses. This normalization step did not affect estimates of the scaling exponents or of the intermittency. Then an *m*-point moving average and an *m*-point moving median were applied to each series of normalized IBIs, resulting in a time series of means, $(\bar{x}_1, \bar{x}_2, ..., \bar{x}_{n+1-m})$, and a time series of medians, $(\tilde{x}_1, \tilde{x}_2, ..., \tilde{x}_{n+1-m})$, for a subject with $n$ IBIs. Fig. 1 illustrates this preprocessing step. For each series of means and medians, $\sigma_k(1/2)$, $\sigma_k(1)$, and $\sigma_k(2)$ were estimated for $k = 1, 2, 4, 8, 16,$ and $32$. Those estimates in turn yielded estimates of $H(1/2)$, $H(1)$, $H(2)$, $\chi(1/2,1)$, and $\chi(1,2)$ using linear regression, as displayed in Fig. 2. The means of estimates of $\sigma_1(1/2)$, $\sigma_1(1)$, $H(1/2)$, and $\chi(1,2)$ are given in Tables 1-2; these four measures were chosen from the others for their superior ability to distinguish heart and stroke subjects, as will be seen.

Thirteen of the measures were used to classify each subject as negative (cardiac) or positive (stroke): $\sigma_1(1/2)$, $\sigma_1(1)$, $\sigma_1(2)$, $\sigma_4(1/2)$, $\sigma_4(1)$, $\sigma_4(2)$, $\sigma_{32}(1/2)$, $\sigma_{32}(2)$, $H(1/2)$, $H(1)$, $H(2)$, $\chi(1/2,1)$, and $\chi(1,2)$. (These measures were selected such that $q=1/2$, $q=1$, and $q=2$ since $\sigma_k(2)$ is a standard deviation of increments, since $H(1)$ is related to a fractal dimension and is a measure of nonstationarity [9], since $H(2)$ is related to the exponent of the power spectrum [9], and since $\chi(1/2,1)$ and $\chi(1,2)$ perform well in simulations [18]. I used $\sigma_1(q)$ because it is computed from the first difference of the time series and I used $\sigma_4(q)$ and $\sigma_{32}(q)$ because second-order dispersions at scales of 4 intervals [19] and 32 intervals [3] effectively discriminate between heart failure and healthy subjects for long time series. If an estimate of $\sigma_k(q)$ was less than or equal to a threshold, a positive (stroke) classification was made, but if the estimate was greater than the threshold, a negative (cardiac) classification was made. On the other hand, if an estimate of $H(q)$ or $\chi(q_1,q_2)$ was greater than or equal to a threshold, a positive classification was made, but if the estimate was less than the threshold, a negative classification was made. Changing the threshold changes both the true-positive rate (the proportion of stroke subjects correctly classified as stroke subjects) and the false-positive rate (the proportion of cardiac subjects incorrectly classified as stroke subjects), so that the true-positive rate can be plotted as a function of the false-positive rate for each estimator. The area under this curve (AUC) is a nonparametric measure of how well a measure distinguishes two groups [20] and provides a reliable way to judge time-series estimators [3,18]; AUC is 0 if everything is classified incorrectly, 1 if everything is classified correctly, and 1/2 if the estimator classifies as well as chance. Table 3 and Fig. 3 report the AUCs for both estimators of intermittency and for the dispersion and scaling estimators with the highest AUCs.

Fig. 3 indicates that the best measures of dispersion are $\sigma_1(1/2)$ and $\sigma_1(1)$, except in the case of no smoothing $(m=1)$, when $\sigma_4(1/2)$ is best, and that the best scaling exponent is $H(1/2)$. The two measures of intermittency, $\chi(1/2,1)$ and $\chi(1,2)$, have comparable performance, but $\chi(1,2)$ is preferred since it has been studied by simulation and applied to heart rate and stock market data [6]. The moving average gives the intermittency measures much greater performance than does the moving median, but the particular filter used has little effect on the performance of the dispersions and scaling exponents. The poorer performance of intermittency with the moving median may be due to an artificial inflation of intermittency estimates, as Tables 1-2 suggest. For a given value of $m$, Tables 1-2 reveal that while the standard errors of the mean of $\chi(1,2)$ estimates are about equal for the moving average and moving median within each subject group, the differences in the means of $\chi(1,2)$ estimates between the heart and stroke groups are much higher for the moving average than for the moving median. Thus, the higher performance of the moving average results from a greater difference in estimated intermittency between the cardiac and stroke groups, not from a smaller variability in the estimates. It appears that the moving median significantly misrepresents the intermittency inherent in the dynamics, rather than merely minimizing the effects of outlying artifacts that are external to the system. However, the moving median might be appropriate in data sets that are corrupted with extreme measurement errors. Alternatively, this kind of contamination may instead call for methods of outlier rejection, such as the application of an *m*-point moving mode using an outlier-proof estimator of the mode like that of Ref. [13]. That estimator has a lower bias and variance than the median when the proportion of outliers is high [14].

It is evident from Fig. 3a that the application of the moving average substantially improves the performance of the three best measures and that $\sigma_1(1/2)$ and $\chi(1,2)$ perform much better than $H(1/2)$. For the 5-point moving average, the estimates of $\sigma_1(1/2)$ and $\chi(1,2)$ statistically differentiate between stroke and cardiac subjects ($P<3\%$), but the estimates of $H(1/2)$ are not significantly different ($P>20\%$), by the two-sided Wilcoxon rank-sum test. (The *P*-values found are actually lower bounds since they are based on some of the best estimates.) That dispersion, a scale-dependent measure, performs better than scale-free scaling exponents has been shown previously [3], but the superior performance of the scale-free intermittency over the scaling exponent is a new result. This result is consistent with the theory behind finite-difference intermittency estimation, that the difference between two scaling exponents has information not present in either scaling exponent considered by itself [4]. As with longer time series, the group with a higher intermittency has a lower dispersion [4]. The scaling exponent may perform much better for longer time series than for shorter time series since it is much more sensitive to the record length than the dispersion: a study that distinguishes heart-failure patients from healthy subjects indicates that changing the number of interbeat intervals from 100 to 100,000 increases the AUC of a scaling exponent by 0.3, but only increases the AUC of a dispersion estimator by 0.06 [3].

Care must be exercised in interpreting estimates of dispersion, scaling, and intermittency since smoothing introduces a bias in parameter estimation. However, this bias poses no more of a problem to comparative studies than does the bias in the estimation of scaling exponents due to the finite size of observed time series. In each case, bias is handled by taking advantage of the fact that time series that are realizations of the same stochastic process will have equal biases, so that the expectation value of the difference in biased estimates is equal to the difference in the estimated parameters. For example, to determine whether a time series can be considered a random walk, the estimate of $H(2)$ should not be compared to 1/2, the value of $H(2)$ for a random walk, since the estimate is biased and its expectation value is not 1/2. The estimate should instead be compared to estimates of $H(2)$ from simulated random walks of the same length as the time series since such estimates will have the same bias [5]. Likewise, the stroke and cardiac subjects can be compared using biased estimates of the parameters since, if there is no difference in the heart rate dynamics between the two groups, then their biases, being equal, will cancel. That there are notable differences in the dispersion and intermittency estimates between these two groups is evidence that the heart rate dynamics differs for cardiac and stroke subjects.

Based on the results of this study, it is hypothesized that, after applying a 5-point moving average, $\sigma_1(1/2)$ and $\chi(1,2)$ can effectively predict health status for other data sets. (In clinical practice, these measures would be used with other tests to increase the true-positive rate and decrease the false-positive rate.) It is likely that different values of *m* are optimal for data with different amounts of noise, or for samples drawn from other populations. For example, since the equivalent of $H(2)$ is higher for older adults [21] and for preterm neonates [5] than for younger adults, heart rate patterns depend on age and different window sizes may better quantify HRV for different ages. This notwithstanding, Fig. 3a shows that the performances of $\sigma_1(1/2)$ and $\chi(1,2)$ are not very sensitive to the choice of *m*, so that reliable results can be obtained using any of a number of different window widths.

I would like to thank Ryutaro Takahashi of the Tokyo Metropolitan Institute of Gerontology for providing the data analyzed.

| m  | Filter | $\sigma_1(1/2)$ | $\sigma_1(1)$ | $H(1/2)$ | $\chi(1,2)$ |
|----|--------|-----------------|---------------|----------|-------------|
| 1  | None   | 15.356+/-2.22   | 20.861+/-2.963 | 0.258+/-0.028 | 0.065+/-0.024 |
| 5  | MA     | 4.445+/-0.452   | 5.73+/-0.587   | 0.548+/-0.025 | 0.054+/-0.021 |
|    | MM     | 2.844+/-0.603   | 6.135+/-1.207  | 0.723+/-0.035 | 0.184+/-0.016 |
| 9  | MA     | 2.881+/-0.28    | 3.692+/-0.359  | 0.658+/-0.022 | 0.049+/-0.019 |
|    | MM     | 1.563+/-0.39    | 4.017+/-0.827  | 0.909+/-0.036 | 0.258+/-0.017 |
| 13 | MA     | 2.194+/-0.219   | 2.759+/-0.269  | 0.722+/-0.02  | 0.04+/-0.016  |
|    | MM     | 1.024+/-0.284   | 3.+/-0.622     | 1.04+/-0.038  | 0.318+/-0.017 |
| 17 | MA     | 1.801+/-0.181   | 2.244+/-0.225  | 0.765+/-0.02  | 0.038+/-0.014 |
|    | MM     | 0.749+/-0.182   | 2.388+/-0.428  | 1.104+/-0.038 | 0.333+/-0.02  |
| 21 | MA     | 1.523+/-0.153   | 1.888+/-0.188  | 0.792+/-0.019 | 0.034+/-0.014 |
|    | MM     | 0.592+/-0.157   | 1.98+/-0.339   | 1.173+/-0.038 | 0.359+/-0.022 |
| 25 | MA     | 1.314+/-0.132   | 1.622+/-0.163  | 0.814+/-0.018 | 0.032+/-0.013 |
|    | MM     | 0.471+/-0.124   | 1.715+/-0.288  | 1.221+/-0.038 | 0.379+/-0.023 |

**Table 1.** Mean +/- standard error of mean of estimates for cardiac subjects. "MA" is the *m*-point moving average and "MM" is the *m*-point moving median.

| m  | Filter | $\sigma_1(1/2)$ | $\sigma_1(1)$ | $H(1/2)$ | $\chi(1,2)$ |
|----|--------|-----------------|---------------|----------|-------------|
| 1  | None   | 10.605+/-1.074  | 16.391+/-2.093 | 0.277+/-0.029 | 0.098+/-0.026 |
| 5  | MA     | 3.11+/-0.265    | 4.418+/-0.427  | 0.595+/-0.018 | 0.149+/-0.031 |
|    | MM     | 1.681+/-0.165   | 3.893+/-0.309  | 0.753+/-0.025 | 0.205+/-0.023 |
| 9  | MA     | 2.169+/-0.182   | 2.946+/-0.263  | 0.695+/-0.016 | 0.139+/-0.03  |
|    | MM     | 0.882+/-0.112   | 2.49+/-0.212   | 0.955+/-0.022 | 0.262+/-0.017 |
| 13 | MA     | 2.169+/-0.182   | 2.257+/-0.202  | 0.757+/-0.015 | 0.127+/-0.029 |
|    | MM     | 0.882+/-0.112   | 1.917+/-0.167  | 1.074+/-0.021 | 0.312+/-0.016 |
| 17 | MA     | 1.38+/-0.133    | 1.82+/-0.171   | 0.801+/-0.014 | 0.124+/-0.03  |
|    | MM     | 0.408+/-0.065   | 1.562+/-0.153  | 1.175+/-0.023 | 0.36+/-0.016  |
| 21 | MA     | 1.167+/-0.117   | 1.52+/-0.147   | 0.83+/-0.013  | 0.12+/-0.03   |
|    | MM     | 0.313+/-0.056   | 1.313+/-0.138  | 1.247+/-0.031 | 0.39+/-0.022  |
| 25 | MA     | 1.03+/-0.104    | 1.333+/-0.13   | 0.843+/-0.014 | 0.112+/-0.03  |
|    | MM     | 0.253+/-0.045   | 1.166+/-0.119  | 1.276+/-0.028 | 0.404+/-0.021 |

**Table 2.** Mean +/- standard error of mean of estimates for stroke subjects. "MA" is the *m*-point moving average and "MM" is the *m*-point moving median.

| m | Filter | $\sigma_1(1/2)$ | $\sigma_1(1)$ | $\sigma_4(1/2)$ | $H(1/2)$ | $H(1)$ | $\chi(1/2,1)$ | $\chi(1,2)$ |
|---|---|---|---|---|---|---|---|---|
| 1 | None | 0.631 | 0.607 | 0.72 | 0.512 | 0.49 | 0.616 | 0.595 |
| 5 | MA | 0.716 | 0.658 | 0.628 | 0.614 | 0.541 | 0.69 | 0.679 |
|   | MM | 0.679 | 0.688 | 0.694 | 0.533 | 0.493 | 0.597 | 0.552 |
| 9 | MA | 0.673 | 0.647 | 0.607 | 0.592 | 0.516 | 0.671 | 0.69 |
|   | MM | 0.658 | 0.698 | 0.647 | 0.645 | 0.56 | 0.605 | 0.546 |
| 13 | MA | 0.645 | 0.626 | 0.586 | 0.601 | 0.516 | 0.701 | 0.709 |
|   | MM | 0.654 | 0.671 | 0.601 | 0.527 | 0.51 | 0.595 | 0.484 |
| 17 | MA | 0.645 | 0.628 | 0.577 | 0.599 | 0.524 | 0.694 | 0.701 |
|   | MM | 0.677 | 0.679 | 0.616 | 0.667 | 0.59 | 0.673 | 0.643 |
| 21 | MA | 0.65 | 0.631 | 0.578 | 0.607 | 0.552 | 0.694 | 0.707 |
|   | MM | 0.66 | 0.66 | 0.611 | 0.614 | 0.563 | 0.631 | 0.614 |
| 25 | MA | 0.633 | 0.612 | 0.578 | 0.594 | 0.524 | 0.664 | 0.673 |
|   | MM | 0.652 | 0.637 | 0.599 | 0.592 | 0.526 | 0.622 | 0.569 |

**Table 3.** AUC (ability of each estimator to distinguish cardiac subjects from stroke subjects). "MA" is the *m*-point moving average and "MM" is the *m*-point moving median.

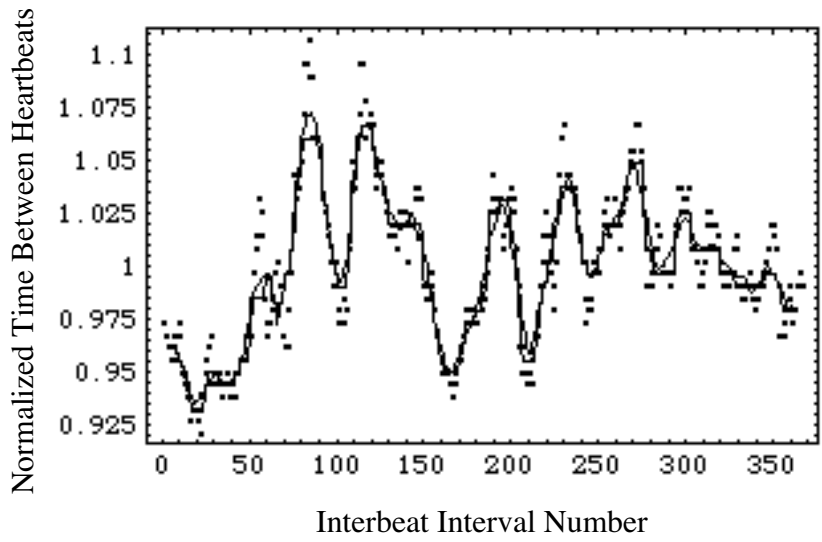

**Fig. 1.** Normalized interbeat-interval time series of a stroke subject. The smooth line is the 13-point moving average and the irregular line is the 13-point moving median of the series. The greatest differences between the two smoothed signals appear near the local extrema.

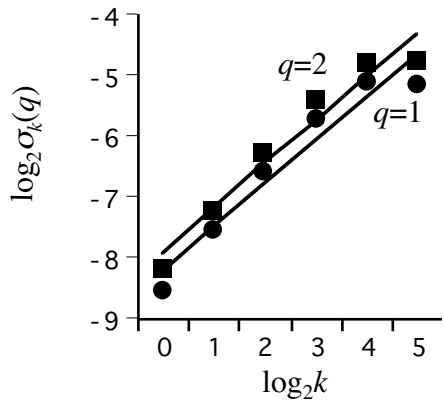

**Fig. 2.** Log-log plot of dispersion versus time scale for the 13-point moving average of the normalized interbeat-interval time series of a stroke subject. Squares represent dispersion estimates for $q=2$ and circles for $q=1$. The slope of the top line is the estimate of $H(2)$ and the slope of the bottom line is the estimate of $H(1)$; the difference between the slopes yields an estimate of intermittency (2).

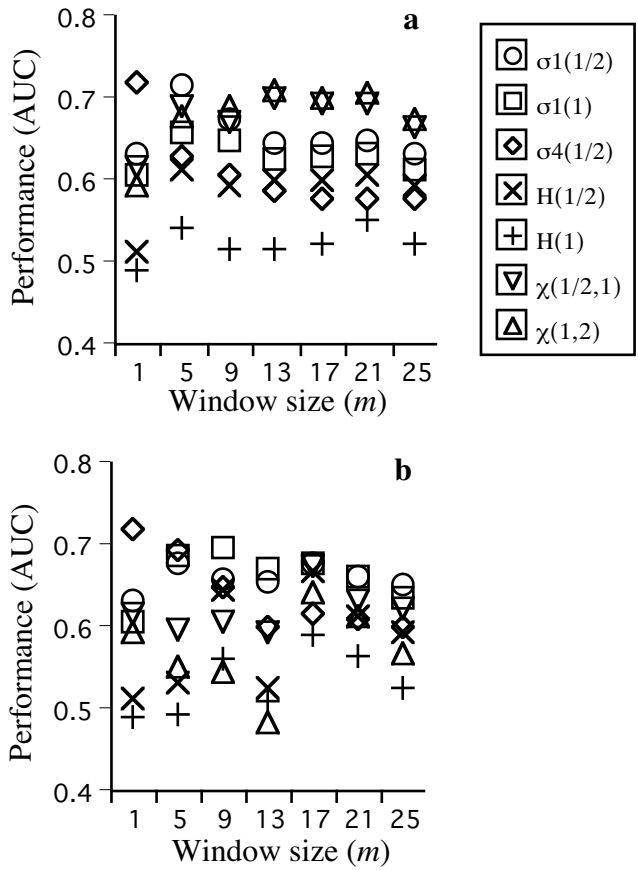

**Fig. 3.** Ability of each estimator to discriminate between cardiac and stroke subjects, with (a) an *m*-point moving average and (b) an *m*-point moving median.

# References

[1] B. J. West and A. L. Goldberger, *American Scientist* **75**, 354 (1987).
[2] P. Ch. Ivanov, L. A. N. Amaral, A. L. Goldberger, S. Havlin, M. G. Rosenblum, Z. Struzik, and H. E. Stanley, *Nature* **399**, 461 (1999).
[3] S. Thurner, M. C. Feurstein, S. B. Lowen, and M. C. Teich, *Phys. Rev. Lett*. **81**, 5688 (1998).
[4] D. R. Bickel, *Phys. Lett. A* **262**, 251 (1999).
[5] D. R. Bickel, M. T. Verklan, and J. Moon, *Physical Review E* **58**, 6440 (1998).
[6] D. R. Bickel and D. Lai, *Computational Statistics and Data Analysis* **37**, 419-431 (2001).
[7] P. J. Bickel and E. L. Lehmann, *Annals of Statistics* **4**, 1139 (1976).
[8] P. Abry and P. Flandrin, In: *Wavelets in Medicine and Biology*, eds. A. Aldroubi and M. Unser, CRC Press, New York, 413 (1996).
[9] A. Davis, A. Marshak, and W. Wiscombe, In: *Wavelets in Geophysics*, edited by E. Foufoula-Georgiou and P. Kumar, 45 (Academic Press, New York, 1994).
[10] D. R. Bickel, *Physica A* **265**, 634-648 (1999)**.**
[11] Y. Wang, *J. R. Statist. Soc. B* **59**, 603-613 (1997).
[12] D. L. Donoho and P. J. Huber, in *Fetschrift for Erich L. Lehmann*, ed. P. J. Bickel, K. A. Doksum, and J. L. Hodges Jr., Wadsworth: Belmont, CA (1983).
[13] D. R. Bickel, "Robust estimators of the mode and skewness of continuous data," *Computational Statistics and Data Analysis* **39**, 153-163 (2002).
[14] D. R. Bickel, *InterStat*, November 2001, http://statjournals.net:2002/interstat/articles/2001/abstracts/n01001.html-ssi (2001).
[15] C. L. do Lago, V. F. Juliano, and C. Kascheres, *Analytica Chimica Acta* **30**, 281-288 (1995).
[16] M. G. Kendall and A. Stuart, *The Advanced Theory of Statistics*, Vol. 2**,** Hafner Publishing Company: New York (1967).
[17] Y. Ashkenazy, P. Ch. Ivanov, S. Havlin, C.-K. Peng, A. L. Goldberger, H. E. Stanley, *Physical Review Letters* **86**, 1900 (2001).
[18] D. R. Bickel, *Chaos, Solitons & Fractals* **13**, 491-497 (2002).
[19] C. Heneghan, S. B. Lowen, and M. C. Teich, *Proc. 1999 IEEE Int. Conf.* (Phoenix, AZ 1999), paper SPTM-8.2 (1999).
[20] M. S. Pepe, *Journal of the American Statistical Association* **95**, 308 (2000).
[21] N. Iyengar, C. K. Peng, R. Morin, A. L. Goldberger, and L. A. Lipsitz, *American Journal of Physiology* **271**, R1078 (1996).